\documentclass[11pt]{amsart}
\usepackage{yhmath}
\usepackage{amssymb}
\usepackage[T1]{fontenc}
\usepackage[frenchb]{babel}
\usepackage{hyperref}
\hypersetup{
    pdftitle={Sous-groupes compacts d'hom\'{e}omorphismes de la sph\`{e}re}, 
    pdfauthor={Boris Kolev},
    pdfsubject={Mathematics}, 
    pdfkeywords={Compact groups of homeomorphisms, Groups acting on specific manifolds}
    }
\newtheorem{thm}{Th\'{e}or\`{e}me}[section]
\newtheorem{cor}[thm]{Corollaire}
\newtheorem{lem}[thm]{Lemme}
\newtheorem{sublem}[thm]{Sous-lemme}

\theoremstyle{definition}
\newtheorem{defn}[thm]{D\'{e}finition}
\theoremstyle{remark}
\newtheorem{rem}[thm]{Remarque}
\numberwithin{equation}{section}

\newcommand{\abs}[1]{\left\vert#1\right\vert}
\newcommand{\set}[1]{\left\{#1\right\}}
\newcommand{\interior}[1]{\mathrm{int}(#1)}
\newcommand{\Real}{\mathbb R}

\newcommand{\OO}{\mathrm{O}}
\newcommand{\SO}{\mathrm{SO}}
\newcommand{\homeo}{\mathrm{Hom\acute{e}o}}

\begin{document}

\title[Sous-groupes compacts d'hom\'{e}omorphismes de la sph\`{e}re]{Note sur les sous-groupes compacts d'hom\'{e}omorphismes de la sph\`{e}re}
\author{Boris Kolev}
\address{CMI,39, rue F. Joliot-Curie, 13453 Marseille cedex 13, France}
\email{boris.kolev@cmi.univ-mrs.fr}

\thanks{Je tiens \`{a} remercier le rapporteur pour sa relecture extr\^{e}mement
attentive de cet article et ses nombreuses remarques qui m'ont aid\'{e} \`{a} am\'{e}liorer ce texte.}%

\subjclass{57s10 (57s25)}

\begin{abstract}
L'objet de cet article est d'exposer la d\'{e}monstration du fait que
tout sous-groupe compact d'hom\'{e}omorphismes de la sph\`{e}re est
topologiquement conjugu\'{e} \`{a} un sous-groupe ferm\'{e} du groupe orthogonal
$\OO(3)$.
\end{abstract}

\maketitle


\section{Introduction}
Le r\'{e}sultat que nous nous proposons d'exposer dans cet article, \`{a}
savoir que tout sous-groupe compact d'hom\'{e}omorphismes de $S^2$ est
topologiquement conjugu\'{e} \`{a} un sous-groupe ferm\'{e} du groupe orthogonal
$\OO(3)$ se situe dans le cadre plus g\'{e}n\'{e}ral d'une suite de
questions connue sous le nom de \emph{5\ieme{} probl\`{e}me de Hilbert}
\cite{Ser52,Yan76}. Plus pr\'{e}cis\'{e}ment, soit $G$ un groupe localement
compact qui agit fid\`{e}lement sur une vari\'{e}t\'{e} $M$, on se pose les
questions suivantes:
\begin{enumerate}
  \item $G$ est-il n\'{e}cessairement localement euclidien\footnote{~Une autre terminologie pour d\'{e}signer une vari\'{e}t\'{e} topologique.} ?
  \item Si $G$ est localement euclidien, est-ce un groupe de Lie ?
  \item Si $G$ est un groupe de Lie, existe-il une structure
  analytique sur $M$ invariante par $G$ ?
\end{enumerate}
La r\'{e}ponse \`{a} la premi\`{e}re question n'est pas connue en dehors de
quelques cas particuliers. La r\'{e}ponse \`{a} la deuxi\`{e}me question est
positive (cf. Gleason \cite{Gle50}, Montgomery and Zippin
\cite{MZ55}). La r\'{e}ponse \`{a} la troisi\`{e}me question est n\'{e}gative en
g\'{e}n\'{e}ral. Il existe des contre-exemples simples dans le cas non
compact. Citons \'{e}galement la construction par Bing \cite{Bin52}
d'une involution n\'{e}gative de $S^3$ non conjugu\'{e}e \`{a} un \'{e}l\'{e}ment de
$\OO(4)$, d'exemples d'hom\'{e}omorphisme p\'{e}riodique de $S^3$ non
conjugu\'{e} \`{a} un \'{e}l\'{e}ment de $\SO(4)$ (Bing \cite{Bin64},
Montgomery-Zippin \cite{MZ55}, Bredon \cite{Bre68}), d'une action de
$\mathbb{U}(1)$ sur $S^4$ non lin\'{e}arisable (Montgomery, Zippin
\cite{MZ55}) et d'une action de $\mathbb{U}(1)$ sur $S^{2n+2}$ non
lin\'{e}arisable \cite{CG97}. Signalons enfin une preuve par Cairns et
Ghys \cite{CG97} que toute action topologique de $\SO(n)$ sur
$\mathbb{R}^n$ qui pr\'{e}serve l'origine est globalement conjugu\'{e}e \`{a}
l'action standard.

Un groupe de Lie poss\`{e}de une propri\'{e}t\'{e} remarquable: il existe un
voisinage de l'identit\'{e} qui ne contient aucun sous-groupe non
trivial. D'un groupe topologique qui poss\`{e}de cette propri\'{e}t\'{e}, on dit
qu'il \emph{n'a pas de petit sous-groupe}. Imm\'{e}diatement apr\`{e}s la
d\'{e}monstration par Haar, en 1933, d'une mesure invariante sur tout
sous-groupe localement compact, von Neumann \cite{Neu33} \'{e}tablit, en
utilisant la th\'{e}orie des repr\'{e}sentations unitaires, le r\'{e}sultat
suivant (voir \'{e}galement \cite{Pon66}), consid\'{e}r\'{e} comme la premi\`{e}re
\'{e}tape majeure dans la r\'{e}solution du 5\ieme{} probl\`{e}me:

\begin{thm}[von Neumann]\label{thm:vonNeumann}
Un groupe \emph{compact} qui ne \emph{poss\`{e}de pas de petit
sous-groupe} est un groupe de Lie.
\end{thm}

Le but de cet article est de pr\'{e}senter une demonstration compl\`{e}te et
moderne d'un th\'{e}or\`{e}me d\^{u} \`{a} Ker\'{e}kj\'{a}rt\'{o} \cite{Ker41} qui donne une
caract\'{e}risation topologique compl\`{e}te du groupe des rotations et de
ses sous-groupes ferm\'{e}s.

\begin{thm}[Ker\'{e}kj\'{a}rt\'{o}]\label{thm:Kerekjarto}
Tout sous-groupe compact de $\homeo(S^{2})$ est topologiquement
conjugu\'{e} \`{a} un sous-groupe ferm\'{e} de $\OO(3)$.
\end{thm}

La preuve donn\'{e}e par Ker\'{e}kj\'{a}rt\'{o} consiste \`{a} \'{e}tablir d'abord qu'un
sous-groupe compact d'hom\'{e}omorphismes de la sph\`{e}re qui poss\`{e}dent un
point fixe laisse invariant un disque topologique autour de ce
point, et que tout groupe compact d'hom\'{e}omorphismes du disque est
topologiquement conjugu\'{e} \`{a} un sous-groupe ferm\'{e} du groupe des
isom\'{e}tries euclidiennes $\OO(2)$. Ces r\'{e}sultats essentiels se
trouvent en r\'{e}alit\'{e} dans des travaux ant\'{e}rieurs de Ker\'{e}kj\'{a}rt\'{o}
\cite{Ker34} connus pour \^{e}tre extr\^{e}mement confus. C'est pourquoi
nous en reprenons la d\'{e}monstration compl\`{e}te dans un langage moderne
dans les premi\`{e}res sections. Le reste de la preuve est une \'{e}tude
casuistique qui se base sur la nature et le nombre des sous-groupes
qui fixent un point (les \emph{stabilisateurs}). Dans \cite{Ker41},
le cas qui appara\^{\i}t le plus compliqu\'{e} et qui occupe la majeur partie
de l'article est celui o\`{u} le groupe agit transitivement sur la
sph\`{e}re car Ker\'{e}kj\'{a}rt\'{o} reconstruit dans ce cas \og \`{a} la main \fg ~la
g\'{e}om\'{e}trie euclidienne de la sph\`{e}re. Dans notre d\'{e}monstration, qui
utilise le langage moderne de la g\'{e}om\'{e}trie diff\'{e}rentiel, ce cas est,
au contraire, le plus simple.

L'article originale de Ker\'{e}kj\'{a}rt\'{o} \cite{Ker41} traite \'{e}galement des
sous-groupes compact d'hom\'{e}omorphismes des autres surfaces
compactes, bien que la majeure partie de l'article soit consacr\'{e}e \`{a}
la sph\`{e}re. En effet, pour les autres surfaces, on se ram\`{e}ne \`{a} des
arguments \'{e}labor\'{e}s pour la sph\`{e}re et le disque. Ainsi, l'\'{e}tude d'un
groupe compact $G$ de transformations d'une surface ferm\'{e}e $M$ de
caract\'{e}ristique d'Euler $\chi(M) \leq 0$, consiste \`{a} traiter d'abord
le sous-groupe $G_{0}$ des transformations isotopes \`{a} l'identit\'{e}. Ce
sous-groupe $G_{0}$ est ferm\'{e}, distingu\'{e} dans $G$ et d'indice fini.
Dans le cas o\`{u} $\chi(M)<0$, on montre en passant au rev\^{e}tement
universel que $G_{0}$ est trivial. Dans le cas du tore, une analyse
analogue \`{a} celle de la section~\ref{sec:circle}, o\`{u} la notion de
\emph{nombre de rotation} est remplac\'{e}e par celle de \emph{vecteur
de rotation}, permet de conclure que $G_{0}$ est conjugu\'{e} \`{a} un
sous-groupe de translation. Cette \'{e}tude n'est pas d\'{e}taill\'{e}e dans cet
article o\`{u} nous nous concentrons essentiellement sur la sph\`{e}re et le
disque. On pourra consulter \cite{BK98} pour plus de d\'{e}tails sur les
autres surfaces.

Un prolongement naturel de ces questions consiste \`{a} rechercher
\'{e}galement une caract\'{e}risation topologique du groupe homographique ou
d'un de ses sous-groupes, question \'{e}galement envisag\'{e}e par
Ker\'{e}kj\'{a}rt\'{o}~\cite{Ker34}. L'\'{e}tude de ce probl\`{e}me a fait appara\^{\i}tre la
notion de groupe de convergence (Ghering and Martin \cite{GM87}).
Mais la r\'{e}ponse ne semble pas aussi simple que pour le groupe des
rotations.

La section~\ref{sec:metric_spaces} de cet article est consacr\'{e}e \`{a}
quelques propri\'{e}t\'{e}s g\'{e}n\'{e}rales des sous-groupes compacts
d'hom\'{e}omorphismes d'un espace m\'{e}trique $(X,d)$ et \`{a} l'\'{e}tude locale
(au voisinage d'un point fixe) lorsque $X$ est une surface. Dans la
section~\ref{sec:circle}, on d\'{e}taille la classification compl\`{e}te des
sous-groupes compacts d'hom\'{e}omorphismes du cercle et dans la
section~\ref{sec:disk}, celle des sous-groupes compacts
d'hom\'{e}omorphismes du disque. La section~\ref{sec:Newman} pr\'{e}sente
une d\'{e}monstration \'{e}l\'{e}mentaire (d\^{u} \`{a} M. H. A. Newman) du fait qu'un
sous-groupe compact d'hom\'{e}omorphismes de la sph\`{e}re n'a pas de petit
sous-groupe. Enfin, la section~\ref{sec:proof} contient l'\'{e}tude
d\'{e}taill\'{e}e des sous-groupes compacts d'hom\'{e}omorphismes de la sph\`{e}re.


\section{Sous-groupes compacts d'hom\'{e}omorphismes d'un espace m\'{e}trique compact}
\label{sec:metric_spaces}

Soit $(X,d)$ un espace m\'{e}trique compact. On d\'{e}finit la distance de
deux applications continues $f,g : X \to X$ par la formule :
\begin{equation*}
  d(f,g)= \max_{x\in X}d(f(x),g(x)).
\end{equation*}
Cette distance munit le groupe des hom\'{e}omorphismes de $(X,d)$ d'une
structure de \emph{groupe topologique}. Nous pouvons \'{e}noncer le
r\'{e}sultat suivant:
\begin{thm}\label{thm:Ascoli}
Soit $G$ un sous-groupe ferm\'{e} d'hom\'{e}omorphismes de $(X,d)$. Les
propositions suivantes sont \'{e}quivalentes:
\begin{enumerate}
  \item $G$ est compact.
  \item L'ensemble des \'{e}l\'{e}ments de $G$ forme une famille
  \'{e}quicontinue.
  \item Il existe une distance sur $X$ pour laquelle les \'{e}l\'{e}ments
  de $G$ sont des isom\'{e}tries.
\end{enumerate}
\end{thm}

\begin{proof}
$(1)\Rightarrow (3)$ est une cons\'{e}quence de l'existence de la mesure
de Haar sur $G$. En effet, ceci nous permet de construire une
distance invariante en prenant la \og moyenne \fg ~ pour la mesure
de Haar des images de la distance $d$ par les \'{e}l\'{e}ments de $G$.
$(3)\Rightarrow (2)$ est trivial et $(2)\Rightarrow (1)$ est un
corollaire direct du th\'{e}or\`{e}me d'Ascoli.
\end{proof}

En particulier, l'ensemble des it\'{e}r\'{e}s d'un \'{e}l\'{e}ment $f$ appartenant \`{a}
un groupe compact d'hom\'{e}omorphismes forme une famille \'{e}quicontinue.
Nous introduirons la d\'{e}finition suivante:

\begin{defn}
Un hom\'{e}omorphisme $f$ d'un espace m\'{e}trique compact $(X,d)$ est
\emph{r\'{e}gulier} si la famille des it\'{e}r\'{e}s de $f$ est \'{e}quicontinue,
autrement dit si $\forall \varepsilon
> 0$, $\exists\alpha
>0$ tel que:
\begin{equation}
d(x,y) < \delta \Rightarrow d(f^n(x),f^ny)) < \varepsilon, \quad
\forall n .
\end{equation}
\end{defn}

Citons quelques exemples: un hom\'{e}omorphisme p\'{e}riodique, une
isom\'{e}trie pour la distance $d$ sont des hom\'{e}omorphismes r\'{e}guliers.
On peut montrer le r\'{e}sultat suivant \cite{BK98} :

\begin{lem}\label{lem:group_closure}
Soit $f$ un hom\'{e}omorphisme r\'{e}gulier, alors la fermeture du groupe
engendr\'{e} par $f$ est compact.
\end{lem}

Dans le cas o\`{u} $X$ est la sph\`{e}re $S^{2}$, ou plus g\'{e}n\'{e}ralement une
surface compacte, nous pouvons expliciter compl\`{e}tement la dynamique
d'un groupe compact d'hom\'{e}omorphisme au voisinage d'un point fixe.

\begin{lem}\label{lem:local_conectivity}
Soit $G$ un sous-groupe compact d'hom\'{e}omorphismes de la sph\`{e}re
$S^{2}$ et $D \subset S^{2}$ un disque topologique ferm\'{e}. Alors le
compact
\begin{equation*}
K = \bigcup_{g\in G}g(D)
\end{equation*}
est localement connexe.
\end{lem}

\begin{proof}
Commen\c{c}ons par rappeler qu'un espace m\'{e}trique compact est localement
connexe si et seulement si pour tout $\varepsilon > 0$, on peut
l'\'{e}crire comme une r\'{e}union finie de compacts connexes de diam\`{e}tre
inf\'{e}rieur \`{a} $\varepsilon$.

Soit $\varepsilon > 0$. Choisissons une triangulation de $D$ dont
les cellules $e_{1},e_{2},\dotsc , e_{r}$, sont de diam\`{e}tre
inf\'{e}rieur \`{a} $\varphi (\varepsilon)$, o\`{u} $\varphi (\varepsilon)$ est
la borne sup\'{e}rieure des nombres positifs $\alpha$ tel que:
\begin{equation*}
    d(x,y)<\alpha\Rightarrow d(g(x),g(y))<\varepsilon,
\end{equation*}
pour $x,y \in S^{2}$ et $g \in G$.

Soit $\rho > 0$, tel que l'int\'{e}rieur de toute cellule $e_{i}$
contienne une boule $B(x_{i}, \rho)$. Alors
\begin{equation}
B(g(x_{i}),\varphi (\rho )) \subset g(e_{i}),\quad \forall i,\quad
\forall g .
\end{equation}
Par cons\'{e}quent, l'aire de chaque cellule $g(e_{i})$ est minor\'{e}e par
$4\pi \sin \varphi (\rho )$ et la famille $\set{ g(e_{i}) }_{i,g}$
ne contient seulement qu'un nombre fini de cellules deux \`{a} deux
disjointes.

Dans cette famille, soit $\set{e_{1}^{\prime}, \dotsc,
e_{p}^{\prime} }$ une sous-famille finie, de cardinal maximal, de
cellules deux \`{a} deux disjointes. Alors pour tout $g \in G$ et tout
$i \in \set{1, \dotsc, r  }$, il existe $j \in \set{1, \dotsc, p }$
tel que $e_{j}^{\prime}\cap g(e_{i})\neq\emptyset$. Pour $k \in
\set{1, \dotsc, p }$, notons $M_{k}$ la fermeture de l'union de
toutes les cellules $g(e_{i})$ qui rencontrent $e_{k}^{\prime}$.
Alors $M_{k}$ est un compact connexe de diam\`{e}tre inf\'{e}rieur \`{a}
$3\varepsilon$ (le diam\`{e}tre de chaque cellule $g(e_{i})$ \'{e}tant
major\'{e} par $\varepsilon$) et
\begin{equation}
    K=\bigcup_{k=1}^{p}M_{k},
\end{equation}
ce qui ach\`{e}ve la d\'{e}monstration.
\end{proof}

\begin{thm}\label{thm:disc_neighbourhoods}
Soit $G$ un sous-groupe compact d'hom\'{e}omorphismes de la sph\`{e}re
$S^{2}$ et $x_{0}$ un point fixe de $G$. Alors il existe un syst\`{e}me
fondamental de voisinage de $x_{0}$, constitu\'{e} par des disques
topologiques invariants par $G$.
\end{thm}

\begin{rem}
D'apr\`{e}s le Th\'{e}or\`{e}me~\ref{thm:Ascoli}, $G$ laisse invariant une
distance $\delta$. On peut donc \^{e}tre tenter de croire que les boules
(pour cette distance $\delta$), centr\'{e}es au point $x_{0}$,
fournissent ce syst\`{e}me de disques invariants. Mais ceci ne
fonctionne pas car on ne sait rien, a priori, de cette distance
invariante $\delta$ obtenue en moyennant par $G$ la distance
euclidienne: elle n'a pas de raison d'\^{e}tre riemannienne.
\end{rem}

La d\'{e}monstration du Th\'{e}or\`{e}me~\ref{thm:disc_neighbourhoods} repose
sur un r\'{e}sultat classique de topologie du plan qui se d\'{e}montre \`{a} la
fois par des m\'{e}thodes purement topologiques \cite{New92,Why64} et
par des m\'{e}thodes issues de la g\'{e}om\'{e}trie complexe \cite{Pom92}.

\begin{thm}\label{thm:Caratheodory}
Soit $K$ un compact, connexe, localement connexe, non vide de
$S^{2}$, non r\'{e}duit \`{a} un point et sans point de coupure\footnote{Un
point $x\in K$ est un \emph{point de coupure} si $K \setminus
\set{x}$ n'est pas connexe.}. Alors la fronti\`{e}re de chaque
composante de $S^{2} \setminus K$ est une courbe ferm\'{e}e simple.
\end{thm}

\begin{proof}[D\'{e}monstration du Th\'{e}or\`{e}me~\ref{thm:disc_neighbourhoods}]
Donnons-nous arbitrairement $\varepsilon >0$. Nous pouvons trouver
$\delta > 0$ tel que:
\begin{equation}\label{equ:first_ball}
    d(x,y) \leq \delta \Rightarrow d(g(x),g(y)) < \varepsilon, \qquad
    \forall x,y \in S^{2}, \; \forall g \in G ,
\end{equation}
puis $\eta > 0$ tel que
\begin{equation}\label{equ:second_ball}
    d(x,y) \leq \eta \Rightarrow d(g(x),g(y)) < \delta, \qquad
    \forall x,y \in S^{2}, \; \forall g \in G .
\end{equation}
Soit $D$ le disque euclidien (ferm\'{e}) de centre $x_{0}$ et de rayon
$\eta$. Formons le compact, connexe, invariant:
\begin{equation*}
  K = \bigcup_{g\in G} g(D).
\end{equation*}
D'apr\`{e}s \eqref{equ:first_ball}, on a $K \subset B(x_{0},\delta)$.
D\'{e}signons par $V_{\infty}$ la composante de $S^{2}\setminus K$ qui
contient $S^{2}\setminus B(x_{0},\delta)$. Soit $g \in G$. En vertu
de \eqref{equ:second_ball}, on a
\begin{equation*}
    g(S^{2}\setminus B(x_{0},\varepsilon))\subset S^{2}\setminus
    B(x_{0},\delta).
\end{equation*}
Par cons\'{e}quent:
\begin{equation*}
    g(S^{2}\setminus B(x_{0},\varepsilon))\subset V_{\infty} \cap
    g(V_{\infty}),
\end{equation*}
et donc $g(V_{\infty}) = V_{\infty}$.

Par ailleurs, d'apr\`{e}s le lemme~\ref{lem:local_conectivity}, $K$ est
localement connexe. On pourra v\'{e}rifier que l'adh\'{e}rence d'un ouvert
connexe, non vide de la sph\`{e}re ne poss\`{e}de pas de point de coupure.
Il en est donc ainsi de $K$ qui est l'adh\'{e}rence de l'ouvert connexe,
non vide
\begin{equation*}
  U = \bigcup_{g\in G} g(\interior{D}).
\end{equation*}
Par suite, chaque composante connexe du compl\'{e}mentaire de $K$ est un
disque topologique en vertu du th\'{e}or\`{e}me~\ref{thm:Caratheodory}. En
particulier, la fronti\`{e}re de $V_{\infty}$ est une courbe ferm\'{e}e
simple invariante et le disque topologique bord\'{e} par cette courbe et
contenant $x_{0}$ est invariant et contenu dans la boule
$B(x_{0},\varepsilon)$.
\end{proof}


\section{Sous-groupes compacts d'hom\'{e}omorphismes du cercle}
\label{sec:circle}

Commen\c{c}ons par d\'{e}montrer les r\'{e}sultats \'{e}l\'{e}mentaires suivants:

\begin{lem}\label{lem:regular_in_1D}
Soit $f:[0,1]\to[0,1]$ un hom\'{e}omorphisme \emph{r\'{e}gulier} et
croissant, alors $f = Id$.

Soit $f:S^{1}\to S^{1}$ un hom\'{e}omorphisme \emph{r\'{e}gulier}, qui
pr\'{e}serve l'orientation et poss\`{e}de un point fixe, alors $f = Id$.
\end{lem}

\begin{proof}
Soit $f:[0,1]\to[0,1]$ un hom\'{e}omorphisme \emph{r\'{e}gulier} et
croissant. Par l'absurde, supposons $f\ne Id$, et soit  $]a,b[$ une
composante de $[0,1]\setminus Fix(f)$. On a $f(a) = a$, $f(b) = b$
et (par exemple):
\begin{equation*}
    f(x) > x , \quad \forall x \in ]a,b[ .
\end{equation*}
Alors, l'orbite par $f$ de tout point de $]a,b[$ converge vers $b$,
ce qui entre en contradiction avec la r\'{e}gularit\'{e} de $f$ qui impose
que l'orbite d'un point voisin de $a$ reste proche de $a$.

Soit maintenant $f:S^{1}\to S^{1}$ un hom\'{e}omorphisme
\emph{r\'{e}gulier}, qui pr\'{e}serve l'orientation et poss\`{e}de un point
fixe. L'\'{e}tude d'un rel\`{e}vement de $f$,
\begin{equation*}
    \tilde{f}:\Real\to\Real,
\end{equation*}
qui poss\`{e}de un point fixe $\tilde{x}$, nous ram\`{e}ne au r\'{e}sultat
pr\'{e}c\'{e}dent en consid\'{e}rant la restriction de $\tilde{f}$ \`{a}
l'intervalle $[\tilde{x},\tilde{x}+1]$.
\end{proof}

Soit $f$ un hom\'{e}omorphisme de $S^{1}$ qui pr\'{e}serve l'orientation et
$\tilde{f}: \mathbb{R} \to \mathbb{R}$ un rel\`{e}vement de $f$. On
rappelle que la limite
\begin{equation}\label{equ:rotation_number_def1}
    \theta(\tilde{f}) = \lim_{n} \; \frac{\tilde{f}^{n}(\tilde{x}) - \tilde{x}}{n}
\end{equation}
existe toujours est ne d\'{e}pend pas du point $\tilde{x}\in \mathbb{R}$
(voir \cite{Dev86}). Si $\tilde{f}^{\prime}$ est un autre rel\`{e}vement
de $f$ alors $\theta(\tilde{f}) - \theta(\tilde{f}^{\prime})$ est un
entier. On note $\rho(f)$ la classe r\'{e}siduelle de ces nombres modulo
$\mathbb{Z}$ et on l'appelle le \emph{nombre de rotation} de $f$.

\begin{lem}\label{lem:rotnumber_morphism}
Soit $G$ un sous-groupe compact d'hom\'{e}omorphismes du cercle $S^{1}$
qui pr\'{e}servent l'orientation, alors l'application \emph{nombre de
rotation}
\begin{equation*}
    \rho : G \to \mathbb{U}(1)
\end{equation*}
est un morphisme continu et injectif.
\end{lem}

\begin{proof}
En moyennant les images par $G$ de la mesure canonique de $S^{1}$ \`{a}
l'aide de la mesure de Haar sur $G$, on obtient une mesure de
probabilit\'{e} $\mu$ sur le cercle invariante par $G$.

Soit $f\in G$. On peut r\'{e}\'{e}crire
l'expression~\eqref{equ:rotation_number_def1} sous la forme
\begin{equation}\label{equ:rotation_number_def2}
    \theta(\tilde{f}) = \lim_{n} \; \sum_{k=0}^{n-1}\frac{\varphi_{\tilde{f}} \, (f^{n}(x))}{n}
\end{equation}
o\`{u} $\varphi_{\tilde{f}}: S^{1} \to \mathbb{R}$ est la fonction
induite par l'application p\'{e}riodique
\begin{equation*}
    \tilde{x} \mapsto \tilde{f}(\tilde{x}) -\tilde{x}, \qquad \tilde{x}\in\Real .
\end{equation*}
Par cons\'{e}quent, d'apr\`{e}s le \emph{th\'{e}or\`{e}me ergodique de Birkhoff}
(voir par exemple \cite{Wal82}), on a
\begin{equation}\label{equ:rotation_number_def3}
    \theta(\tilde{f}) = \int_{S^{1}} \varphi_{\tilde{f}} \; d\mu \;.
\end{equation}

Enfin, on d\'{e}montre sans difficult\'{e}, \`{a} partir de la d\'{e}finition de
$\varphi_{\tilde{f}}$, la \emph{relation de cocycle}:
\begin{equation*}
    \varphi_{\tilde{f}\circ\tilde{g}} = \varphi_{\tilde{f}} \circ g +
    \varphi_{\tilde{g}} .
\end{equation*}
Par suite, si $f$ et $g$ sont deux \'{e}l\'{e}ments quelconques du groupe
$G$, on obtient
\begin{equation}\label{equ:commutativity}
    \theta(\tilde{f} \circ \tilde{g}) = \int_{S^{1}} \varphi_{\tilde{f}} \circ g \; d\mu
    + \int_{S^{1}} \varphi_{\tilde{g}} \; d\mu = \theta(\tilde{f})
    + \theta(\tilde{g}) ,
\end{equation}
ce qui \'{e}tablit que $\rho$ est bien un morphisme de groupe.

L'injectivit\'{e} est une cons\'{e}quence du Lemme~\ref{lem:regular_in_1D}
et du fait qu'un hom\'{e}omorphisme du cercle, qui a pour nombre de
rotation $0$, poss\`{e}de un point fixe (voir \cite{Dev86}).

La continuit\'{e} r\'{e}sulte de l'in\'{e}galit\'{e} suivante:
\begin{equation*}
     d(f,Id) \leq \varphi(\varepsilon) \Rightarrow \abs{ \theta (\tilde{f}) } \leq \varepsilon,
\end{equation*}
qui est elle-m\^{e}me cons\'{e}quence de l'\'{e}quicontinuit\'{e} d'un \'{e}l\'{e}ment $f$
du groupe $G$.
\end{proof}

\section{Sous-groupes compacts d'hom\'{e}omorphismes du disque}
\label{sec:disk}

Le r\'{e}sultat suivant g\'{e}n\'{e}ralise un r\'{e}sultat connu pour les
hom\'{e}omorphismes p\'{e}riodiques du disque \cite{CK94}.

\begin{lem}\label{lem:identity_on_boundary}
Un hom\'{e}omorphisme r\'{e}gulier du disque $D^{2}$ qui est l'identit\'{e} sur
le bord est l'identit\'{e} sur le disque tout entier.
\end{lem}

\begin{proof}
Formons le \emph{double} de $f$, qui est un hom\'{e}omorphisme de la
sph\`{e}re et que nous continuerons de d\'{e}signer par $f$. Nous obtenons
ainsi un hom\'{e}omorphisme r\'{e}gulier, qui est l'identit\'{e} sur une courbe
ferm\'{e}e simple $J$ (correspondant au bord de $D^{2}$) et que nous
pouvons supposer \^{e}tre l'\'{e}quateur.

Choisissons sur $J$ deux points diam\'{e}tralement oppos\'{e}s que nous
noterons $a$ et $b$. Soit $d$ un cercle arbitraire s\'{e}parant les
points $a$ et $b$. Reprenons la construction donn\'{e}e dans la
d\'{e}monstration du Th\'{e}or\`{e}me~\ref{thm:disc_neighbourhoods}, en prenant
pour $G$, la fermeture dans $\homeo(S^{2})$ du groupe engendr\'{e} par
$f$ et pour $D$, le disque (ferm\'{e}) d\'{e}limit\'{e} par $d$ et contenant
$a$. Notons, comme pr\'{e}c\'{e}demment
\begin{equation*}
  K = \bigcup_{g\in G} g(D)
\end{equation*}
et d\'{e}signons par $\Delta$ la composante de $b$ dans $S^{2}\setminus
K$. Alors $\Delta$ est un disque topologique invariant bord\'{e} par une
courbe ferm\'{e}e simple que nous noterons $\delta$. Cette courbe s\'{e}pare
les points $a$ et $b$ et rencontre donc n\'{e}cessairement la courbe
$J$. Par cons\'{e}quent, en utilisant \`{a} nouveau le
lemme~\ref{lem:regular_in_1D}, on en d\'{e}duit que $f=Id$ sur $\delta$.
Or, par construction
\begin{equation*}
    \delta \subset \bigcup_{g \in G}g(d),
\end{equation*}
et donc $\delta \subset d$ (car $\delta\subset Fix(G)$), ce qui
n'est possible que si $\delta=d$. Le cercle $d$ ayant \'{e}t\'{e} choisi
arbitrairement, ceci montre que $f = Id$ sur $S^{2}$.
\end{proof}

\begin{cor}\label{cor:unique_fixed_point}
Soit $f\in \homeo^{+}(D^{2})$ un hom\'{e}omorphisme r\'{e}gulier, diff\'{e}rent
de l'identit\'{e}. Alors $f$ poss\`{e}de un point fixe unique et ce point
est situ\'{e} \`{a} l'int\'{e}rieur du disque.
\end{cor}

\begin{proof}
D'apr\`{e}s le th\'{e}or\`{e}me du point fixe de Brouwer, $f$ poss\`{e}de au moins
un point fixe $x_{0}$. Si $f\ne Id$, il r\'{e}sulte du
lemme~\ref{lem:identity_on_boundary} et du
lemme~\ref{lem:regular_in_1D} que ce point fixe se trouve \`{a}
l'int\'{e}rieur du disque.

Nous allons maintenant montrer que si $f$ poss\`{e}de un second point
fixe $x_{1}$ alors $f=Id$. Pour cela, construisons \`{a} l'aide du
th\'{e}or\`{e}me~\ref{thm:disc_neighbourhoods}, une courbe ferm\'{e}e simple
invariante $J \in \interior{D^{2}}$ qui s\'{e}pare les deux points fixes
et borde un disque topologique (ouvert) $\Delta$ contenant $x_{0}$.
Par construction, l'anneau topologique (ferm\'{e}) $A = D^{2}\setminus
\Delta$ et invariant par $f$ et contient l'autre point fixe,
$x_{1}$. Soit
\begin{equation*}
    \tilde{f}: \Real\times [0,1] \to \Real\times [0,1]
\end{equation*}
un rel\`{e}vement de la restriction de $f$ \`{a} $A$. On peut v\'{e}rifier que
$\tilde{f}$ est un hom\'{e}omorphisme r\'{e}gulier de $\Real\times [0,1]$
(pour la m\'{e}trique standard). Soit $\varphi_{\tilde{f}}: A \to
\mathbb{R}$ la fonction d\'{e}finie par
\begin{equation*}
    x \mapsto p_{1}\left(\tilde{f}(\tilde{x}) - \tilde{x}\right) ,
\end{equation*}
o\`{u} $p_{1}$ est la projection sur le premier facteur du produit
$\Real\times [0,1]$. La r\'{e}gularit\'{e} de $f$ implique l'existence et
l'unicit\'{e} (ind\'{e}pendance par rapport \`{a} $x$) de la limite
\begin{equation*}
    \lim_{n} \; \sum_{k=0}^{n-1}\frac{\varphi_{\tilde{f}} \, ( f^{n}(x) )}{n} ,
\end{equation*}
que nous noterons $\theta(\tilde{f})$ comme dans la preuve du
lemme~\ref{lem:rotnumber_morphism}. Soit $\tilde{x}_{1}$ un
rel\`{e}vement de $x_{1}$ et choisissons pour $\tilde{f}$ un rel\`{e}vement
de $f$ qui fixe $\tilde{x}_{1}$. Alors on a n\'{e}cessairement
$\theta(\tilde{f})=0$ et ceci impose \`{a} $f$ d'avoir un point fixe sur
$\partial D^{2}$ (voir \cite{Dev86}). Il r\'{e}sulte alors du
lemme~\ref{lem:regular_in_1D}, que $f$ est l'identit\'{e} sur $\partial
D^{2}$ puis que que $f$ est l'identit\'{e} sur le disque d'apr\`{e}s le
lemme~\ref{lem:identity_on_boundary}.
\end{proof}

\begin{cor}\label{cor:restriction_morphism}
Soit $G$ un sous-groupe compact d'hom\'{e}omorphismes du disque $D^{2}$.
La restriction au bord
\begin{equation*}
    R : G \rightarrow \homeo(\partial D^{2})
\end{equation*}
est un morphisme continu et injectif.
\end{cor}

\begin{proof}
La restriction au bord d'un sous-groupe d'hom\'{e}omorphismes du disque
est toujours un morphisme continu mais n'est pas injectif en
g\'{e}n\'{e}ral. Soit $g \in G$ un \'{e}l\'{e}ment du noyau de $R$. Comme $g$ est
r\'{e}gulier, alors $g=Id$ en vertu du
lemme~\ref{lem:identity_on_boundary} et donc $R$ est injectif si $G$
est compact.
\end{proof}

En combinant le corollaire~\ref{cor:restriction_morphism} avec le
lemme~\ref{lem:rotnumber_morphism}, on obtient:

\begin{cor}\label{cor:U1}
Tout sous-groupe compact $G$ de $\homeo^{+}(D^{2})$ est isomorphe
(en tant que groupe topologique) \`{a} un sous-groupe ferm\'{e} de
$\mathbb{U}(1)$.
\end{cor}

Si $G$ est fini, il est engendr\'{e} par un \'{e}l\'{e}ment p\'{e}riodique $f$. Dans
ce cas, on montre que $f$ est conjugu\'{e} \`{a} une rotation euclidienne
(voir \cite{CK94}). Sinon, $G = \mathbb{U}(1)$ et nous allons
\'{e}tablir le r\'{e}sultat suivant:

\begin{thm}\label{thm:U1_action}
Toute action continue et fid\`{e}le du groupe $\mathbb{U}(1)$ sur le
disque est topologiquement conjugu\'{e}e \`{a} l'action standard de
$\SO(2)$.
\end{thm}

Nous diviserons la d\'{e}monstration de ce r\'{e}sultat en deux \'{e}tapes: nous
montrerons d'abord que la structure des orbites d'un tel groupe est
identique \`{a} celle de l'action standard et ensuite, ce qui est la
partie la plus d\'{e}licate, qu'il existe un \emph{arc transverse} aux
orbites, ce qui nous permettra de conclure.

\begin{lem}\label{lem:U1_orbits_structure}
Les orbites de toute action continue et fid\`{e}le de $\mathbb{U}(1)$
sur le disque $D^{2}$ sont constitu\'{e}es par un point fixe unique
$x_{0}$ \`{a} l'int\'{e}rieur du disque et des courbes ferm\'{e}es simples qui
entourent ce point.
\end{lem}

\begin{proof}
Soit $G$ l'image de $\mathbb{U}(1)$ dans $\homeo(D^{2})$ et $f\in G$
un \'{e}l\'{e}ment d'ordre infini. L'unique point fixe $x_{0}$ de $f$ donn\'{e}
par le corollaire~\ref{cor:unique_fixed_point} est \'{e}galement un
point fixe de $G$ car les it\'{e}r\'{e}s de $f$ forment un ensemble dense
dans $G$. Par cons\'{e}quent, $x_{0}$ est \'{e}galement l'unique point fixe
de tout autre \'{e}l\'{e}ment $g\ne Id$ de $G$. Il en r\'{e}sulte que la
$G$-orbite de tout point $x$ diff\'{e}rent de $x_{0}$ est une courbe
ferm\'{e}e simple. Cette courbe est invariante par $f$ et borde un
disque qui contient n\'{e}cessairement un point fixe de $f$, en vertu du
th\'{e}or\`{e}me du point fixe de Brouwer. Ce point fixe ne peut \^{e}tre que
$x_{0}$, ce qui ach\`{e}ve la d\'{e}monstration du
lemme~\ref{lem:U1_orbits_structure}.
\end{proof}

\begin{lem}\label{lem:transveral_arc}
\'{E}tant donn\'{e} une action topologique et fid\`{e}le de $\mathbb{U}(1)$ sur
le disque $D^{2}$, il existe un arc simple rencontrant chaque orbite
en un point et un seul.
\end{lem}

\begin{rem}
Je ne connais pas de preuve \'{e}l\'{e}mentaire de ce lemme. La
d\'{e}monstration propos\'{e}e ici est une construction \og \`{a} la main \fg
~de cet \emph{arc transverse}. On pourra remarquer que ce r\'{e}sultat
n'est pas une cons\'{e}quence imm\'{e}diate du
lemme~\ref{lem:U1_orbits_structure}. Il existe en effet des
partitions de l'anneau par des familles de courbes ferm\'{e}es simples
essentielles qui n'admettent pas de transversale. On pourra
consulter \cite{Whi41} pour plus de d\'{e}tails sur le sujet.
\end{rem}

Soit $G$ l'image de $\mathbb{U}(1)$ dans $\homeo(D^{2})$ et $x_{0}$,
l'unique point fixe de $G$. Pour tout $x \in D^{2}$, on note
$\gamma(x)$ la $G$-orbite de $x$.

\begin{sublem}\label{ass:small_arc}
Pour tout $\varepsilon > 0$, il existe $\delta > 0$ tel que si $x$
et $y$ sont deux points distincts d'une m\^{e}me $G$-orbite $\gamma$ et
$d(x, y) < \delta$, alors un des deux arcs d\'{e}limit\'{e}s par $x$ et $y$
sur $\gamma$ a un diam\`{e}tre inf\'{e}rieur \`{a} $\varepsilon$.
\end{sublem}

\begin{proof}
Notons $x_{0}$, l'unique point fixe de $G$. Soit $\varepsilon > 0$
et $\Delta$ un disque contenant $x_{0}$, invariant par $G$ et de
diam\`{e}tre inf\'{e}rieur \`{a} $\varepsilon$ (voir
th\'{e}or\`{e}me~\ref{thm:disc_neighbourhoods}). On pose $A = D^{2}\setminus
\Delta$. Il suffit de d\'{e}montrer le sous-lemme pour les orbites
contenues dans l'anneau $A$, ce qui r\'{e}sulte de l'observation
suivante. Il existe un voisinage ouvert connexe $V$ de l'identit\'{e}
dans $G$ tel que
\begin{equation}
d(x,g(x)) < \varepsilon /2 ,
\end{equation}
pour tout $x\in A$ et $g\in V$. Comme de plus $G$ agit librement sur
$A$, il existe $\delta > 0$ tel que:
\begin{equation}
d(x,g(x)) \geq \delta ,
\end{equation}
pour tout $x \in A$ et $g \in G \setminus V$. Par cons\'{e}quent, si $x$
et $y$ sont deux points d'une m\^{e}me $G$-orbite $\gamma\subset A$ tels
que $d(x, y) ) < \delta$ alors $y=g(x)$ avec $g\in V$ et l'arc
$\set{g(x); \, g\in V}$ de $\gamma$ a un diam\`{e}tre inf\'{e}rieur \`{a}
$\varepsilon$.
\end{proof}

On munit l'ensemble des orbites de $G$ d'un ordre total en
d\'{e}finissant la relation suivante: $\gamma \leq \gamma^{'}$
(respectivement $\gamma < \gamma^{'}$) si et seulement si $\gamma$
est contenue dans le disque ferm\'{e} (respectivement ouvert) bord\'{e} par
$\gamma^{'}$ (avec la convention $x_{0} \leq x_{0}$).

\begin{defn}
\'{E}tant donn\'{e} deux points $x,y \in D^{2}$, tel que $\gamma(x) <
\gamma(y)$, on appelle \textit{$\mu$-cha\^{\i}ne monotone} de $x$ \`{a} $y$,
une collection $\set{x_{0}=x, x_{1},\dotsc, x_{n}=y}$ de points de
$D^{2}$ tels que:
\begin{equation*}
    d(x_{k}, x_{k+l}) < \mu \quad \text{et} \quad \gamma(x_{k}) <
\gamma (x_{k+l}),
\end{equation*}
pour $k=0, \dotsc , n-1$.
\end{defn}

\begin{sublem}\label{ass:monotone_chain}
Pour tout $\varepsilon > 0$, il existe $\delta > 0$ tel que deux
points quelconques $x,y$ n'appartenant pas \`{a} la m\^{e}me orbite et
v\'{e}rifiant $d(x, y) < \delta$ peuvent \^{e}tre joints par une
$\mu$-cha\^{\i}ne monotone de diam\`{e}tre inf\'{e}rieur \`{a} $\varepsilon$ pour
tout $\mu > 0$.
\end{sublem}

\begin{proof}
Soit $\varepsilon > 0$ et choisissons $\delta > 0$ ($\delta <
\varepsilon$) comme dans le sous-lemme~\ref{ass:small_arc}. Soient
$x,y$ deux points tels que $\gamma(x) < \gamma(y)$ et $d(x, y) <
\delta/2$. \'{E}tant donn\'{e} $\mu > 0$ ($\mu < \delta$), on peut trouver
une suite finie de $G$-orbites
\begin{equation*}
\gamma_{0} = \gamma (x) < \gamma_{1} < \dotsb < \gamma_{n}
=\gamma(y)
\end{equation*}
telle que la \emph{distance de Hausdorff}\footnote{La distance de
Hausdorff est d\'{e}finie sur les parties ferm\'{e}es d'un espace m\'{e}trique
compact $(X,d)$ par la formule:
\begin{equation*}
    d_{H}(A, B) = \max \set{ \max_{a\in A} d (a, B), \,  \max_{b\in B} d (b, A) }.
\end{equation*}
} $d_{H}(\gamma_{k}, \gamma_{k+1})$ de deux courbes cons\'{e}cutives de
la suite soit inf\'{e}rieure \`{a} $\mu/3$.

Le segment $xy$ rencontre chacune des courbes interm\'{e}diaires
$\gamma_{k}$. Choisissons pour chaque $k\in\set{1,\dotsc ,n-1}$, un
point
\begin{equation*}
    x_{k} \in xy \cap\gamma_{k} .
\end{equation*}
Si $d(x_{k}, x_{k+1}) < \mu$ pour tout $k$, nous avons construit
notre $\mu$-cha\^{\i}ne de diam\`{e}tre inf\'{e}rieur \`{a} $\varepsilon$, sinon,
voil\`{a} comment raffiner cette cha\^{\i}ne pour en obtenir une.

Soit $x_{r},x_{r+1}$, une paire de points cons\'{e}cutifs tels que
$d(x_{r}, x_{r+1}) \geq \mu$. Comme $d_{H}(\gamma_{r}, \gamma_{r+1})
< \mu/3$, on peut trouver un point $x^{\prime}_{r+1}$ sur
$\gamma_{r+1}$ tel que:
\begin{equation*}
    d(x_{r}, x^{'}_{r+1}) < \mu/3 ,
\end{equation*}
et donc:
\begin{equation*}
    d(x_{r+1}, x^{'}_{r+1}) \leq d(x_{r+1}, x_{r}) +  d(x_{r}, x^{'}_{r+1}) <
    \delta/2 + \mu/3 < \delta .
\end{equation*}
Alors, d'apr\`{e}s le sous-lemme~\ref{ass:small_arc}, un des deux arcs
d\'{e}limit\'{e}s par $x_{r+1}$ et $x^{'}_{r+1}$ sur $\gamma_{r+1}$ a un
diam\`{e}tre inf\'{e}rieur \`{a} $\varepsilon$. Subdivisons cet arc en $s$
sous-arcs de diam\`{e}tre inf\'{e}rieur \`{a} $\mu/3$ et notons les point
interm\'{e}diaires de la fa\c{c}on suivante:
\begin{equation*}
x^{'}_{r+1} = z^{0}_{r+1},z^{1}_{r+1}, \dotsc ,z^{s}_{r+1}= x_{r+1}.
\end{equation*}
Choisissons ensuite des courbes
\begin{equation*}
    \gamma_{r} =  \gamma^{0} <  \gamma^{1} < \dotsb <  \gamma^{s} =
    \gamma_{r+1} .
\end{equation*}
Comme $d_{H}(\gamma_{r}, \gamma_{r+1}) < \mu/3$ et $\gamma_{r} <
\gamma^{j} < \gamma_{r+1}$ pour $1 \leq j \leq s-1$, il est possible
de trouver, sur chaque courbe interm\'{e}diaire $\gamma^{j}$, un point
$x^{j}_{r+1}$ tel que $d(x^{j}_{r+1}, z^{j}_{r+1}) < \mu/3$. La
suite
\begin{equation*}
    x_{r} = x^{0}_{r+1}, x^{1}_{r+1},\dotsc, z^{s}_{r+1} = x_{r+1}
\end{equation*}
est donc une $\mu$-cha\^{\i}ne joignant $x_{r}$ et $x_{r+1}$ qui
appartient \`{a} un $2\varepsilon$-voisinage du segment $x_{r}x_{r+1}$.
En effectuant les corrections n\'{e}cessaires pour chaque paire
$x_{r},x_{r+1}$ telle que $d(x_{r}, x_{r+1}) \geq \mu$, on obtient
finalement une $\mu$-cha\^{\i}ne monotone joignant $x$ et $y$ et de
diam\`{e}tre inf\'{e}rieur \`{a} $4\varepsilon$.
\end{proof}

\begin{proof}[Preuve du Lemme~\ref{lem:transveral_arc}]
En utilisant le sous-lemme~\ref{ass:monotone_chain}, on peut
construire une suite de nombres r\'{e}els $\delta_{n} > 0$, qui tend
vers $0$, et telle que deux points quelconques $x,y$ v\'{e}rifiant $d(x,
y) < \delta_{n}$, peuvent \^{e}tre joints, pour tout $\mu > 0$, par une
$\mu$-cha\^{\i}ne de diam\`{e}tre inf\'{e}rieur \`{a} $1/2^{n}$.

Soit $X_{0}$ une $\delta_{0}$-cha\^{\i}ne monotone joignant $x_{0}$, le
point fixe du groupe \`{a} un point $x_{\infty}$ sur le bord du disque
$D^{2}$. R\'{e}cursivement, ayant d\'{e}fini $X_{n}$, on joint chaque paire
de points cons\'{e}cutifs ${x^{n}_{k}, x^{n}_{k+1}}$ de $X_{n}$, par une
$\delta_{n+1}$-cha\^{\i}ne monotone de diam\`{e}tre inf\'{e}rieur \`{a} $1/2^{n}$
afin d'obtenir $X_{n+1}$ et on pose:
\begin{equation*}
X = \bigcup_{n\in\mathbb{N}}X_{n}.
\end{equation*}
C'est alors un exercice standard de topologie (\cite[Theorem
2.27]{HY88}) de montrer que l'adh\'{e}rence $\overline{X}$ de $X$ dans
$D^{2}$ est un arc joignant $x_{0}$ et $x_{\infty}$. Cet arc est
simple et rencontre chaque orbite en un point unique, par
construction, ce qui ach\`{e}ve la d\'{e}monstration.
\end{proof}

\begin{proof}[Preuve du th\'{e}or\`{e}me~\ref{thm:U1_action}]
Pour compl\'{e}ter la preuve du th\'{e}or\`{e}me~\ref{thm:U1_action}, on choisit
un arc $\alpha$, transverse aux orbites, donn\'{e} par le
lemme~\ref{lem:transveral_arc} et une param\'{e}trisation $x(r), r \in
[0,1]$ de cet arc. L'application
\begin{equation*}
h : re^{i\theta} \mapsto \Psi(e^{i\theta}, x (r)),
\end{equation*}
o\`{u} $\Psi : \mathbb{U}(1)\times D^{2} \to D^{2}$ d\'{e}note l'action,
nous donne alors une conjugaison topologique avec le groupe des
rotations euclidiennes, $\SO(2)$.
\end{proof}


\section{Un lemme de M.H.A. Newman}
\label{sec:Newman}

Avant d'entreprendre l'\'{e}tude des sous-groupes compacts
d'hom\'{e}omorphismes de la sph\`{e}re, nous pr\'{e}sentons un lemme d\^{u} \`{a} M.H.A.
Newman \cite{New31}.
\begin{lem}\label{lem:Newman}
Soit $f$ un hom\'{e}omorphisme p\'{e}riodique de $S^2$ de p\'{e}riode $p>1$,
alors parmi les it\'{e}r\'{e}s de $f$, il en existe un, disons $f^{r}$, tel
que :
\begin{equation}\label{equ:Newman1}
d(f^{r},Id) > 1 .
\end{equation}
De plus :
\begin{equation}\label{equ:Newman2}
d(f,Id) > \frac{2}{p}.
\end{equation}
\end{lem}

\begin{proof}
Commen\c{c}ons par remarquer que \eqref{equ:Newman2} est une cons\'{e}quence
de \eqref{equ:Newman1}. En effet, dans \eqref{equ:Newman1}, on peut
supposer $r\leq p/2$ car l'in\'{e}galit\'{e}~\eqref{equ:Newman1} est
\'{e}quivalente \`{a}
\begin{equation*}
d(f^{p-r},Id) > 1 \; .
\end{equation*}
Par cons\'{e}quent la n\'{e}gation de \eqref{equ:Newman2} conduit \`{a}
\begin{equation*}
d(f^{r}(x),x) \leq \sum_{i=0}^{r-1} d(f^{i+1}(x),f^{i}(x)) \leq
\frac{2r}{p} \leq 1 , \quad \forall x ,
\end{equation*}
et donc \`{a} la n\'{e}gation de \eqref{equ:Newman1}.

La preuve de la premi\`{e}re in\'{e}galit\'{e} r\'{e}sulte de la remarque suivante :
supposons au contraire que $d(f^k,Id) \leq 1$ pour tout k, alors
l'orbite d'un point quelconque $x$ est enti\`{e}rement contenue dans
l'h\'{e}misph\`{e}re de p\^{o}le $x$ et par suite, pour tout p-uplet
$(\lambda_0,\lambda_1,\dotsc,\lambda_{p-1})$ de nombres positifs
tels que $\Sigma \lambda_i = 1$:
\begin{equation}
g_{\lambda}(x) = \sum_{i=0}^{p-1}\lambda_{i}f^{i}(x)\neq 0 ,
\end{equation}
pour tout $x$. Ceci implique l'existence d'une homotopie dans
$\mathbb{R}^3-\{ 0 \}$ entre l'identit\'{e} et la fonction
\begin{equation}
  g(x) = \frac{1}{p} \sum_{i=0}^{p-1} f^i(x),
\end{equation}
ce qui est incompatible avec le fait que $\deg (g)= 0$ modulo $p$.
\end{proof}

\begin{cor}
La boule unit\'{e} ferm\'{e}e de $\homeo(S^{2})$ ne contient aucun
sous-groupe compact non trivial.
\end{cor}

\begin{proof}
En effet, si un tel groupe existe, on peut trouver un \'{e}l\'{e}ment non
trivial $f$ de ce groupe tel que
\begin{equation*}
    d(f^{n},Id) \leq 1 ,
\end{equation*}
pour tout $n\in\mathbb{Z}$. En vertu du lemme~\ref{lem:Newman}, $f$
n'est pas p\'{e}riodique. Quitte \`{a} remplacer $f$ par son carr\'{e}, on peut
supposer que $f$ pr\'{e}serve l'orientation et donc poss\`{e}de au moins un
point fixe. Alors, d'apr\`{e}s le lemme~\ref{thm:disc_neighbourhoods},
$f$ laisse invariant un disque topologique et la fermeture du groupe
engendr\'{e} par $f$ est isomorphe \`{a} $\mathbb{U}(1)$. Mais ce groupe
contient des \'{e}l\'{e}ment p\'{e}riodiques $g$ v\'{e}rifiant \'{e}galement
\begin{equation*}
    d(g^{n},Id) \leq 1 ,
\end{equation*}
pour tout $n\in\mathbb{Z}$, ce qui est en contradiction avec le
lemme~\ref{lem:Newman}.
\end{proof}

Il en r\'{e}sulte qu'un sous-groupe compact d'hom\'{e}omorphismes de la
sph\`{e}re ne \emph{poss\`{e}de pas de petit sous-groupe}. En vertu du
Th\'{e}or\`{e}me~\ref{thm:vonNeumann}, on peut donc \'{e}noncer:

\begin{thm}\label{thm:Hilbert}
Tout sous-groupe compact d'hom\'{e}omorphismes de la sph\`{e}re est un
groupe de Lie
\end{thm}


\section{Preuve du th\'{e}or\`{e}me principal}
\label{sec:proof}

Cette section est consacr\'{e} \`{a} la d\'{e}monstration du
th\'{e}or\`{e}me~\ref{thm:Kerekjarto}. Nous envisagerons dans un premier
temps le cas d'un sous-groupe compact $G$ qui ne contient que des
\'{e}l\'{e}ments qui pr\'{e}servent l'orientation, puis le cas g\'{e}n\'{e}ral.

\subsection{$G$ ne contient que des \'{e}l\'{e}ments qui pr\'{e}servent
l'orientation}

\begin{lem}\label{lem:rotation}
Soit $G$ un sous-groupe compact de $\homeo^{+}(S^{2})$. Alors tout
\'{e}l\'{e}ment de $G$ est topologiquement conjugu\'{e} \`{a} une rotation
euclidienne d'ordre finie ou non.
\end{lem}

\begin{proof}
Soit $g\in G$ un \'{e}l\'{e}ment non trivial. En tant qu'hom\'{e}omorphisme qui
pr\'{e}serve l'orientation de la sph\`{e}re, $g$ poss\`{e}de un point fixe
$x_{0}$ (th\'{e}or\`{e}me de Lefschetz, par exemple). Du
th\'{e}or\`{e}me~\ref{thm:disc_neighbourhoods}, on d\'{e}duit l'existence d'un
disque invariant $\Delta$ contenant $x_{0}$. Le disque
$\overline{S^{2}\setminus\Delta}$ est \'{e}galement invariant et
contient donc un second point fixe $x_{0}^{*}$ de $g$. En vertu du
corollaire~\ref{cor:unique_fixed_point}, $Fix(g)$ est r\'{e}duit \`{a} ces
deux points.

Si $g$ est d'ordre fini, il est conjugu\'{e} \`{a} une rotation euclidienne
d'ordre fini (voir \cite{CK94}). Sinon, l'adh\'{e}rence $H$ du groupe
engendr\'{e} par $g$ est isomorphe \`{a} $\mathbb{U}(1)$ (voir
corollaire~\ref{cor:U1}). Les orbites de $H$ sont constitu\'{e}es par
les deux points fixes $x_{0},x_{0}^{*}$ et des courbes ferm\'{e}es
simples qui s\'{e}parent $x_{0}$ et $x_{0}^{*}$. Le
lemme~\ref{lem:transveral_arc} nous assure l'existence d'un arc
transverse aux orbites, joignant $x_{0}$ et $x_{0}^{*}$, ce qui nous
permet, comme dans la preuve du th\'{e}or\`{e}me~\ref{thm:U1_action},
d'\'{e}tablir que $g$ est topologiquement conjugu\'{e} \`{a} une rotation
euclidienne d'ordre infini.
\end{proof}

Soit $x\in S^{2}$. Le \emph{stabilisateur} de $x$, not\'{e} $Stab(x)$,
est le sous-groupe des \'{e}l\'{e}ments $g$ de $G$ tels que $g(x) = x$. Le
sous-groupe compact $Stab(x)$ laisse invariant un disque contenant
$x$. Il est isomorphe \`{a} un sous-groupe ferm\'{e} du groupe
$\mathbb{U}(1)$ en vertu du corollaire~\ref{cor:U1}.

Inversement, soit $H$ un sous-groupe ferm\'{e} de $G$. Si $Fix(H)\ne
\emptyset$ alors $Fix(H)$ contient exactement deux points
$x_{0},x_{0}^{*}$ et $H$ est isomorphe \`{a} un sous-groupe ferm\'{e} du
groupe $\mathbb{U}(1)$.

\begin{lem}
Si $G$ est infini, il poss\`{e}de un sous-groupe isomorphe \`{a}
$\mathbb{U}(1)$.
\end{lem}

\begin{proof}
D'apr\`{e}s le th\'{e}or\`{e}me~\ref{thm:Hilbert}, $G$ est un groupe de Lie. Par
cons\'{e}quent, si $G$ est infini, sa dimension est sup\'{e}rieure \`{a} $1$
(car $G$ est compact). Il poss\`{e}de donc un sous-groupe \`{a} un param\`{e}tre
non trivial et contient des \'{e}l\'{e}ments d'ordre infini. Soit $g$ un tel
\'{e}l\'{e}ment. Alors l'adh\'{e}rence du groupe engendr\'{e} par $g$ est isomorphe
\`{a} $\mathbb{U}(1)$.
\end{proof}

Nous allons maintenant envisager les divers cas possibles.

\subsubsection{Cas~1: G est fini}

Chaque \'{e}l\'{e}ment non trivial de $G$ poss\`{e}de exactement deux points
fixes. Seulement un nombre fini de points de la sph\`{e}re ont un
stabilisateur non trivial. Soit $\Sigma$ cet ensemble, alors la
projection canonique $\pi : S^2 \to S^2/G$ est un rev\^{e}tement ramifi\'{e}
et on a la \emph{formule de Riemann-Hurwitz}:
\begin{equation}\label{equ:Hurwitz}
\chi (S^2) = n \chi(S^2/G) - \sum_{\bar{s}\in \Sigma / G}
(\nu_{\bar{s}}-1)
\end{equation}
o\`{u} $n$ d\'{e}signe le cardinal de $G$ et $\nu_{\bar{s}}$ est le cardinal
des stabilisateurs des points de ramification $s\in\bar{s}$. De
cette formule, il r\'{e}sulte que $\chi(S^2/G) = 2$ et donc que $S^2/G$
est hom\'{e}omorphe \`{a} $S^2$. Ces rev\^{e}tements sont enti\`{e}rement classifi\'{e}s
par l'action du groupe $G$ sur l'ensemble fini $\Sigma$. A chaque
solution donn\'{e}e par la formule~\eqref{equ:Hurwitz} correspond un
sous-groupe fini de $\SO(3)$ et donc une conjugaison topologique de
$G$ avec ce sous-groupe.

\subsubsection{Cas~2: il n'y a qu'un stabilisateur infini}

Soit $H$ ce stabilisateur et d\'{e}signons par $x_{0}$ et $x_{0}^*$ les
points fixes de $H$. Soit $g$ un \'{e}l\'{e}ment quelconque du groupe $G$.
Alors le point $g(x_{0})$ est \'{e}galement d'indice infini et donc
n\'{e}cessairement
\begin{equation*}
    g(x_{0})\in Fix(H) = \set{x_{0}, x_{0}^*} .
\end{equation*}
Si de plus $g(x_{0}) = x_{0}$, alors $g\in H$.

- Si ceci se produit pour tous les \'{e}l\'{e}ments du groupe alors $G=H$ et
le groupe $G$ est topologiquement conjugu\'{e} au sous-groupe des
rotations euclidiennes autour d'un axe donn\'{e}.

- Sinon, on peut trouver un \'{e}l\'{e}ment $\sigma$ dans $G$ tel que
$\sigma(x_{0}) = x_{0}^*$ et $\sigma(x_{0}^*) = x_{0}$. On a alors
n\'{e}cessairement $\sigma^2 = Id$ et $\sigma h \, \sigma = h^{-1}$,
pour tout $ h \in H$. $G$ est donc isomorphe au groupe di\'{e}dral
infini
\begin{equation*}
D_{\infty} =\mathbb{Z}_{2} \ltimes \mathbb{U}(1)
\end{equation*}
De plus, $\sigma$ induit une involution continue sur l'espace des
orbites de $H$, qui est hom\'{e}omorphe \`{a} un intervalle, et \'{e}change ses
deux extr\'{e}mit\'{e}s. Par suite, $\sigma$ laisse invariant une et une
seule des orbites de $H$. Les points fixes de $\sigma$ sont
n\'{e}cessairement sur cette courbe. On construit alors facilement une
conjugaison entre ce groupe $G$ est l'action euclidienne standard du
groupe $D_{\infty}$.

\subsubsection{Cas~3: il y a au moins deux stabilisateurs infinis distincts}

Dans ce cas, le groupe $G$ agit transitivement sur $S^2$ en
vertu du r\'{e}sultat suivant.

\begin{lem}
S'il existe deux stabilisateurs infinis distincts (en tant que
sous-groupes de $G$), alors le groupe $G$ agit transitivement sur la
sph\`{e}re.
\end{lem}

\begin{proof}
Soit $a$ et $b$ deux points d'indice infini ayant des stabilisateurs
distincts $H_{a}$ et $H_{b}$ respectivement. L'orbite du point $a$
sous l'action du groupe $H_{b}$, que l'on notera $H_b(a)$ est une
courbe ferm\'{e}e simple passant par $a$ et qui rencontre toutes les
orbites du sous-groupe $H_a$ assez voisines de $a$. Par suite, la
$G$-orbite du point $a$, $G(a)$, est ouverte et ferm\'{e}e dans $S^{2}$,
ce qui ach\`{e}ve la d\'{e}monstration.
\end{proof}

Fixons donc un point $a$ de la sph\`{e}re et d\'{e}signons par $H$ le
stabilisateur de ce point (n\'{e}cessairement isomorphe \`{a}
$\mathbb{U}(1)$). L'espace homog\`{e}ne, $G/H \cong S^2$, est muni
naturellement d'une structure de vari\'{e}t\'{e} analytique sur laquelle $G$
agit \'{e}galement de fa\c{c}on analytique. On a donc trouv\'{e} sur $S^2$ une
structure analytique invariante par $G$. On peut alors construire,
par moyennisation, une m\'{e}trique riemannienne sur $S^2$ invariante
par $G$. Comme l'action de $G$ est transitive, cette m\'{e}trique est \`{a}
courbure constante. Quitte \`{a} multiplier cette m\'{e}trique par une
constante, on peut supposer que cette courbure est $1$. Par
cons\'{e}quent, cette vari\'{e}t\'{e} riemannienne est isom\'{e}trique \`{a} la sph\`{e}re
standard et cette isom\'{e}trie d\'{e}finit la conjugaison recherch\'{e}e entre
$G$ et $\SO(3)$.

\subsection{$G$ contient des \'{e}l\'{e}ments qui renversent l'orientation}

Commen\c{c}ons par rappeler les faits suivants. Soit $s$ un
hom\'{e}omorphisme r\'{e}gulier de la sph\`{e}re qui renverse l'orientation.

\begin{itemize}
  \item Si $s$ poss\`{e}de un point fixe, alors n\'{e}cessairement $s^{2} = Id$ et
$s$ est topologiquement conjugu\'{e} \`{a} une r\'{e}flexion orthogonale
\cite{CK94}.
  \item Si $s^{2} = Id$ mais $s$ est sans point fixe, alors
$S^{2}/s$ est hom\'{e}omorphe au plan projectif et par suite, du fait de
l'unicit\'{e} du rev\^{e}tement universel, $s$ est topologiquement conjugu\'{e}
\`{a} la sym\'{e}trie centrale $x\mapsto -x$.
\end{itemize}

Soit $G_{0}$ le sous-groupe de $G$ des \'{e}l\'{e}ments qui pr\'{e}servent
l'orientation. $G_{0}$ est un sous-groupe distingu\'{e} d'indice $2$.
Quitte \`{a} effectuer une premi\`{e}re conjugaison, on peut supposer que
$G_{0}\subset \SO(3)$.

\subsubsection{Cas 1: $G_{0} \simeq \SO(3)$}

Soit $s\in G\setminus G_{0}$. Tout cercle (euclidien) de $S^{2}$ est
une orbite du stabilisateur $Stab(x)$ d'un point dans $G_{0}$. Par
cons\'{e}quent, l'image par $s$ de tout cercle de $S^{2}$ est un cercle.
D'apr\`{e}s un r\'{e}sultat bien connu, ceci entra\^{\i}ne que $s$ est une
\emph{anti-homographie} de la sph\`{e}re. Comme de plus $s$ est
r\'{e}gulier, $s$ appartient n\'{e}cessairement \`{a} $\OO(3)$ et ceci permet de
conclure que $G=\OO(3)$

\subsubsection{Cas 2: $G_{0} \simeq \SO(2)$}

D\'{e}signons par $x_{0}$ et $x_{0}^{*}$ les points fixes de $G_{0}$ et
soit $s\in G\setminus G_{0}$. Alors, $s$ permute les orbites de
$G_{0}$ et induit un hom\'{e}omorphisme sur le quotient $S^{2}/G_{0}$
qui est un intervalle. Il y a donc deux possibilit\'{e}s:

- Si $s$ fixe $x_{0}$ et $x_{0}^{*}$, alors $s^{2} = Id$ et $s$ est
conjugu\'{e} \`{a} une r\'{e}flexion. $G$ est un produit semi-direct de
$\mathbb{Z}_{2}$ par $\mathbb{U}(1)$ et $sgs = g^{-1}$ pour tout
$g\in G_{0}$. En effet, sinon on aurait $sgs = g$ pour tout $g\in
G_{0}$ et la courbe $Fix(s)$ qui contient $x_{0}$ et $x_{0}^{*}$
serait une orbite de $G_{0}$ ce qui n'est pas possible. Par
ailleurs, $s$ induit un hom\'{e}omorphisme croissant sur l'intervalle
$S^{2}/G_{0}$ de p\'{e}riode $2$ qui ne peut donc \^{e}tre que l'identit\'{e}.
Par suite, $s$ pr\'{e}serve les orbites de $G_{0}$. Chaque orbite non
triviale de $G_{0}$ rencontre $Fix(s)$ en deux points au moins. Mais
si $x$ et $y$ sont deux points distincts de $Fix(s)$ appartenant \`{a}
la m\^{e}me orbite de $G_{0}$, alors $y = g(x)$ pour un certain $g\in
G_{0}$ et la relation $sgs =g^{-1}$ nous donne $g^{2}(x) = x$. Donc
$g$ est n\'{e}cessairement d'ordre $2$ et ceci nous permet de conclure
que chaque orbite non triviale de $G_{0}$ ne rencontre $Fix(s)$
qu'en deux points seulement. On est alors en mesure de construire
une conjugaison entre $G$ et le sous-groupe de $\OO(3)$ engendr\'{e} par
les rotations autour de l'axe $x_{0}x_{0}^{*}$ et la r\'{e}flexion par
rapport \`{a} un \'{e}quateur contenant cet axe.

- Si $s$ \'{e}change $x_{0}$ et $x_{0}^{*}$, alors $s$ induit un
hom\'{e}omorphisme d\'{e}croissant sur l'intervalle $S^{2}/G_{0}$. Cet
hom\'{e}omorphisme a un unique point fixe qui correspond \`{a} une orbite
$J$ de $G_{0}$ invariante par $s$. Alors, quitte \`{a} composer $s$ avec
une rotation de $G_{0}$ on peut supposer que $s$ a un point fixe sur
$J$ et donc se ramener encore une fois au cas o\`{u} $s^{2}= Id$ et $s$
est conjugu\'{e} \`{a} une r\'{e}flexion. On a $Fix(s) = J$ et donc $sgs = g$
pour tout $g\in G_{0}$. Dans ce cas, $G$ est le produit direct de
$\mathbb{Z}_{2}$ par $\mathbb{U}(1)$ et on construit facilement une
conjugaison entre $G$ et le sous-groupe de $\OO(3)$ engendr\'{e} par les
rotations autour de l'axe $x_{0}x_{0}^{*}$ et la r\'{e}flexion par
rapport \`{a} l'\'{e}quateur orthogonal \`{a} cet axe.

\subsubsection{Cas 3: $G_{0} \simeq  D_{\infty}$}

Dans ce cas $G_{0}$ est engendr\'{e} par le groupe des rotations autour
d'un axe $x_{0}x_{0}^{*}$ et par un retournement $\rho$ qui \'{e}change
$x_{0}$ et $x_{0}^{*}$ et dont l'axe est perpendiculaire \`{a} la droite
$x_{0}x_{0}^{*}$. Soit $s\in G\setminus G_{0}$, alors $s$ permute
\'{e}galement les points $x_{0}$ et $x_{0}^{*}$ car le sous-groupe des
rotations autour de l'axe $x_{0}x_{0}^{*}$ est invariant par $s$.
Quitte \`{a} composer $s$ avec $\rho$, on peut supposer que $s(x_{0}) =
x_{0}$ et $s(x_{0}^{*})=x_{0}^{*}$ et donc que $s$ est une r\'{e}flexion
topologique. En conjuguant $s$ par une rotation d'axe
$x_{0}x_{0}^{*}$, on peut supposer \'{e}galement que $Fix(s)$ contient
les points fixes de $\rho$. Alors $s\rho s = \rho$. On construit
alors facilement une conjugaison entre $G$ et le sous-groupe de
$\OO(3)$ engendr\'{e} par $G_{0}$ et la r\'{e}flexion plane par rapport au
plan contenant l'axe $x_{0}x_{0}^{*}$ et l'axe du retournement
$\rho$.

\subsubsection{Cas 4: $G_{0}$ est fini}

Dans ce cas $G_{0}$ appartient \`{a} un des cinq types bien connus de
sous-groupe fini de $\SO(3)$ \cite{Ber77}. Il y a deux possibilit\'{e}s:

- $G\setminus G_{0}$ ne contient aucune \emph{r\'{e}flexion
topologique}, autrement dit, $Fix(s) = \emptyset$ pour tout $s\in
G\setminus G_{0}$. Alors $S^{2}/G$ est hom\'{e}omorphe au plan projectif
et la projection canonique $\pi : S^2 \to S^2/G$ est un rev\^{e}tement
ramifi\'{e}, ce qui permet de conclure que $G$ est conjugu\'{e} \`{a} un
sous-groupe fini de $\OO(3)$.

- $G\setminus G_{0}$ contient une \emph{r\'{e}flexion topologique}. Dans
ce cas, $G$ est un produit semi-direct
\begin{equation*}
    \mathbb{Z}_{2} \ltimes G_{0}.
\end{equation*}
Par ailleurs, si $s$ et $s^{\prime}$ sont deux r\'{e}flexions
topologiques distinctes dans $G$ alors le cardinal de l'ensemble
$Fix(s)\cap Fix(s^{\prime})$ est \'{e}gal \`{a} $2$. En effet, un point fixe
commun \`{a} $s$ et $s^{\prime}$ est un point fixe de la rotation
$ss^{\prime}$ qui en poss\`{e}de au plus $2$ et si les deux courbes
$Fix(s)$ et $Fix(s^{\prime})$ ne s'intersectent pas, ou seulement en
un point, alors la rotation $ss^{\prime}$ envoie un disque ferm\'{e} \`{a}
l'int\'{e}rieur de lui-m\^{e}me (\`{a} l'exclusion \'{e}ventuellement d'un point du
bord), ce qui n'est pas possible pour une rotation. On est alors en
mesure de construire \og \`{a} la main \fg, dans chacune des cinq
situations possibles, des domaines fondamentaux et de montrer,
chaque fois, que $G$ est conjugu\'{e} \`{a} un sous-groupe fini de $\OO(3)$
(voir \cite{CK94} et \cite{Ker19}).

\bibliographystyle{amsplain}
\bibliography{kerek}

\providecommand{\bysame}{\leavevmode\hbox to3em{\hrulefill}\thinspace}
\providecommand{\MR}{\relax\ifhmode\unskip\space\fi MR }
\providecommand{\MRhref}[2]{%
  \href{http://www.ams.org/mathscinet-getitem?mr=#1}{#2}
}
\providecommand{\href}[2]{#2}
\begin{thebibliography}{10}

\bibitem{Ber77}
M~Berger, \emph{G{\'e}om{\'e}trie}, 2 ed., vol.~1, Cedic/{F}ernand {N}athan,
  Paris, 1979.

\bibitem{Bin52}
R.~H. Bing, \emph{A homeomorphism between the {$3$}-sphere and the sum of two
  solid horned spheres}, Ann. of Math. (2) \textbf{56} (1952), 354--362.
  \MR{14,192d}

\bibitem{Bin64}
\bysame, \emph{Inequivalent families of periodic homeomorphisms of
  {$E\sp{3}$}}, Ann. of Math. (2) \textbf{80} (1964), 78--93. \MR{29 \#611}

\bibitem{BK98}
C.~Bonatti and B.~Kolev, \emph{Surface homeomorphisms with zero dimensional
  singular set}, Topology and its Applications \textbf{90} (1998), 69--95.

\bibitem{Bre68}
G.~E. Bredon, \emph{Exotic actions on spheres}, Proc. Conf. on Transformation
  Groups (New Orleans, La., 1967), Springer, New York, 1968, pp.~47--76.
  \MR{MR0266239 (42 \#1146)}

\bibitem{CG97}
G.~Cairns and E.~Ghys, \emph{The local linearization problem for smooth {${\rm,
  SL}(n)$}-actions}, Enseign. Math. (2) \textbf{43} (1997), no.~1-2, 133--171.
  \MR{98i:57067}

\bibitem{CK94}
A.~Constantin and B.~Kolev, \emph{The theorem of {K}er{\'e}kj{\'a}rt{\'o} on
  periodic homeomorphisms of the disk and the sphere}, L'enseignement
  Math{\'e}matique \textbf{40} (1994), 193--204.

\bibitem{Ker19}
B.~de~Ker{\'e}kj{\'a}rt{\'o}, \emph{{{\"U}}ber die endlichen topologischen
  {G}ruppen der {K}ugelfl{\"a}che}, Proc. Acad. Amsterdam \textbf{22} (1919).

\bibitem{Ker34}
\bysame, \emph{Topologische charakterisierung der linearen {A}bbildungen}, Acta
  scient. math., Szeged \textbf{6} (1934), 235--262, Erg{\"a}nzug. ibid., 7
  (1934) 58--59.

\bibitem{Ker41}
\bysame, \emph{Sur les groupes compacts de transformations topologiques des
  surfaces.}, Acta Math. \textbf{74} (1941), 129--173 (French).

\bibitem{Dev86}
R.L. Devaney, \emph{An introduction to chaotic dynamical systems}, The
  Benjamin/Cummings Publishing Co. Inc., Menlo Park, CA, 1986. \MR{MR811850
  (87e:58142)}

\bibitem{GM87}
F.~W. Gehring and G.~J. Martin, \emph{Discrete quasiconformal groups. {I}},
  Proc. London Math. Soc. (3) \textbf{55} (1987), no.~2, 331--358.
  \MR{88m:30057}

\bibitem{Gle50}
A.~M. Gleason, \emph{Spaces with a compact {L}ie group of transformations},
  Proc. Amer. Math. Soc. \textbf{1} (1950), 35--43. \MR{11,497e}

\bibitem{HY88}
J.G. Hocking and G.S. Young, \emph{Topology}, Dover Publications Inc., New
  York, 1988. \MR{90h:54001}

\bibitem{MZ55}
D.~Montgomery and L.~Zippin, \emph{Topological transformation groups}, Robert
  E. Krieger Publishing Co., Huntington, N.Y., 1974, Reprint of the 1955
  original. \MR{52 \#644}

\bibitem{New31}
M.~H.~A. Newman, \emph{A theorem on periodic transformations of spaces}, Quart.
  J. Math. \textbf{2} (1931), 1--8.

\bibitem{New92}
M.H.A. Newman, \emph{Elements of the topology of plane sets of points}, second
  ed., Dover Publications Inc., New York, 1992. \MR{93d:54002}

\bibitem{Pom92}
Ch. Pommerenke, \emph{Boundary behaviour of conformal maps}, Grundlehren der
  Mathematischen Wissenschaften [Fundamental Principles of Mathematical
  Sciences], vol. 299, Springer-Verlag, Berlin, 1992. \MR{95b:30008}

\bibitem{Pon66}
L.~S. Pontryagin, \emph{Topological groups}, Translated from the second Russian
  edition by Arlen Brown, Gordon and Breach Science Publishers, Inc., New York,
  1966. \MR{34 \#1439}

\bibitem{Ser52}
J.~P. Serre, \emph{Le cinqui{\`e}me probl{\`e}me de {H}ilbert. {E}tat de la
  question en 1951}, Bull. Soc. Math. France \textbf{80} (1952), 1--10.
  \MR{MR0049205 (14,136a)}

\bibitem{Neu33}
J.~von Neumann, \emph{Die {E}inf{\"u}hrung analytischer {P}arameter in
  topologischen {G}ruppen}, Ann. of Math. (2) \textbf{34} (1933), no.~1,
  170--190. \MR{MR1503104}

\bibitem{Wal82}
P.~Walters, \emph{An introduction to ergodic theory}, Graduate Texts in
  Mathematics, vol.~79, Springer-Verlag, New York, 1982. \MR{MR648108
  (84e:28017)}

\bibitem{Whi41}
H.~Whitney, \emph{On regular families of curves}, Bull. Amer. Math. Soc.
  \textbf{47} (1941), 145--147. \MR{2,322b}

\bibitem{Why64}
G.T. Whyburn, \emph{Topological analysis}, Second, revised edition. Princeton
  Mathematical Series, No. 23, Princeton University Press, Princeton, N.J.,
  1964. \MR{29 \#2758}

\bibitem{Yan76}
C.~T. Yang, \emph{Hilbert's fifth problem and related problems on
  transformation groups}, Mathematical developments arising from Hilbert
  problems (Proc. Sympos. Pure Math., Northern Illinois Univ., De Kalb, Ill.,
  1974), Amer. Math. Soc., Providence, R. I., 1976, pp.~142--146. Proc. Sympos.
  Pure Math., Vol. XXVIII. \MR{54 \#13948}

\end{thebibliography}
\end{document}